\documentclass{amsart}
\usepackage{amssymb}
\begin{document}
\title[A combinatorial problem in infinite groups]{A combinatorial problem in infinite groups}%
\author{Alireza Abdollahi}%
\address{Department of Mathematics,University of
Isfahan,Isfahan 81746-73441, Iran.}%
\email{a.abdollahi@sci.ui.ac.ir}%

\begin{abstract}
Let $w$ be a word in the free group of rank $n \in \mathbb{N}$ and
let  $\mathcal{V}(w)$ be the variety of groups defined by the law
$w=1$. Define $\mathcal{V}(w^*)$ to be the class of all groups
$G$ in which for any infinite subsets $X_1, \dots, X_n$ there
exist $x_i \in X_i$, $1\leq i\leq n$, such that $w(x_1, \dots,
x_n)=1$. Clearly, $\mathcal{V}(w) \cup \mathcal{F} \subseteq
\mathcal{V}(w^*)$; $\mathcal{F}$ being the class of finite
groups. In this paper, we investigate some words $w$ and some
certain classes $\mathcal{P}$ of groups for which the equality
$\left(\mathcal{V}(w) \cup \mathcal{F}\right)\cap \mathcal{P}=
\mathcal{P} \cap \mathcal{V}(w^*)$ holds.
\end{abstract}
\maketitle

\begin{center} {\bf\Large Introduction and results} \end{center}

Let $w$ be a word in the free group of rank $n \in \mathbb{N}$ and
let  $\mathcal{V}(w)$ be the variety of groups defined by the law
$w=w(x_1,\dots,x_n)=1$. P. Longobardi, M. Maj and A. Rhemtulla in
\cite{LMR} defined $\mathcal{V}(w^*)$ to be the class of all
groups $G$ in which for any infinite subsets $X_1, \dots, X_n$
there exist $x_i \in X_i$, $1\leq i\leq n$, such that $w(x_1,
\dots, x_n)=1$ and raised the question of whether $\mathcal{V}(w)
\cup \mathcal{F} =\mathcal{V}(w^*)$ is true; $\mathcal{F}$ being
the class of finite groups. There is no example, so far, of an
infinite group in $\mathcal{V}(w^*)\backslash\mathcal{V}(w)$. In
fact the origin of this problem is the following observation:\\
Let $G$ be an infinite group such that in every two infinite
subsets of $G$ there exist two commuting  elements, then $G$ is
abelian. This is an immediate consequence
 of the answer of B. H. Neumman
to a question of P. Erd\"os; B. H. Neumman proved that an infinite
group $G$ is centre-by-finite if and only if every infinite subset of $G$
contains two distinct commuting elements \cite{N}. Since this first paper,
problems of a  similar nature have been the object of several articles (for
example  \cite{A2}, \cite{A3}, \cite{AT1},
 \cite{D1}, \cite{D3}, \cite{DRS}, \cite{G},
\cite{LW}, \cite{LM2}, \cite{LMMR}, \cite{RH}).\\
As far as we
 know,  the equality $\mathcal{V}(w)\cup \mathcal{F}=\mathcal{V}(w^*)$ is known for the
following words:
$w=x^m$, $w=[x_1, \dots, x_n]$ \cite{LMR}, $w=[x,y]^2$ \cite{LM},
$w=[x,y,y]$ \cite{S1}, $w=[x,y,y,y]$ \cite{S2},
$w=(xy)^{-3}x^3y^3$ \cite{A1}, $w=x_1^{\alpha_1}\cdots x_m^{\alpha_m}$
 where $\alpha_1, \dots, \alpha_m$ are non-zero integers \cite{AT2},
$w=(xy)^2(yx)^{-2}$ or  $w=[x^m,y]$ where $m\in\{3,6\}\cup \{2^k
\;|\; k\in\mathbb{N}\}$ \cite{AT3}, $w=[x^n,y][x,y^n]^{-1}$ where
$n\in\{\pm 2,3\}$ \cite{Taeri} and $w=[x^m,y^m]$ or
$w=(x_1^mx_2^m\cdots x_n^m)^2$ where
$m\in\{2^k \;|\; k\in \mathbb{N}\}$ \cite{Bouk}. \\
\par
In \cite{PS}, P. Puglisi and L. S. Spiezia proved that every   infinite
locally finite group
 (or locally soluble group) in $\mathcal{V}([x,_ky]^*)$  is a $k$-Engel
group; (recall that  $[x,_ky]$ is defined inductively by $[x,
_0y]=x$ and $[x,_ky]=[[x,_{k-1}y],y]$ for $k\in\mathbb{N}$). In
\cite{D2},  C. Delizia proved the equality $\mathcal{V}(w)\cup
\mathcal{F}=\mathcal{V}(w^*)$ on the  classes of hyperabelian,
locally soluble and locally finite groups where $w=[x_1, \dots,
x_k, x_1]$ and $k$ is an integer greater than 2. Later G.
Endimioni generalized these results by proving that every
infinite locally finite or locally soluble group in
$\mathcal{V}(w^*)$
 belongs to the variety $\mathcal{V}(w)$, where
$w$ is a word in a free group such that finitely generated
soluble groups in $\mathcal{V}(w)$ are nilpotent (see Theorem 3
of \cite{E}) (recall that the variety $\mathcal{V}([x_1, \dots,
x_k, x_1])$  ($k>2$)
 is exactly the variety of  nilpotent groups of nilpotency class
at most $k$ \cite{Mac} and   every finitely generated soluble
Engel group  is nilpotent \cite{Gru}.)\\
We say that a group $G$ is locally graded
if and only if every finitely generated non-trivial subgroup of $G$ has a
non-trivial finite quotient.
We proved in  Theorem 4 of \cite{A3}
 that an infinite locally graded group in
$\mathcal{V}([x_1,_kx_2]^*)$ is a $k$-Engel group. We generalize
this result as Theorem A, below. In order to state our first
result we need the following definition.
 Following  \cite{KR} we say
 that a group $G$ is restrained  if and only if
 $\left<x\right>^{\left<y\right>}=\left<x^{y^i}
 \;|\; i\in \mathbb{Z}\right>$ is finitely
generated for all $x,y\in G$. We show by Proposition 1 below, why the
following theorem improves the above mentioned results.\\
\par
{\bf Theorem A.} {\sl Let $w$ be a word in a free group such that
every finitely generated residually finite  group in
$\mathcal{V}(w)$ is polycyclic-by-finite. Then every infinite
finitely generated locally graded restrained group
 in $\mathcal{V}(w^*)$  belongs to the variety $\mathcal{V}(w)$.}\\
\par
G. Endimioni proved that every infinite locally nilpotent group in
$\mathcal{V}(w^*)$  belongs to the variety $\mathcal{V}(w)$, where
$w$ is a word in a free group (see Theorem 1 of \cite{E}).
The following theorem generalizes Theorem 1 of \cite{E}.\\
\par
{\bf Theorem B.} {\sl Let $w$ be a word in a free group and let
$\mathcal{P}$ be a  class of groups which satisfies the
following conditions:\\
(1)\; the class  $\mathcal{P}$ is closed under taking subgroups.\\
(2)\; every $\mathcal{P}$-group is  soluble.\\
(3)\; every infinite finitely generated
($\mathcal{P}$-by-finite)-group in
$\mathcal{V}(w^*)$ belongs to the variety $\mathcal{V}(w)$.\\
Then every infinite residually [(locally
$\mathcal{P}$)-by-finite] group in $\mathcal{V}(w^*)$ belongs to
 $\mathcal{V}(w).$\\ }
\par
For example, the classes of nilpotent groups,
 polycyclic groups, abelian-by-nilpotent groups and
  soluble residually finite groups satisfy the assumptions of Theorem B.\\
\par
Here we also  obtain some reductions in investigation of the
equality $\mathcal{V}(w) \cup \mathcal{F}=\mathcal{V}(w^*)$ on
certain classes of groups and
certain words $w$. For example\\
\par
{\bf Theorem C.} {\sl Let $w$ be a  non-trivial word
 in a free group. Then every non-linear simple locally finite
 group does not  belong to the class $\mathcal{V}(w^*)$.}\\
\par
In \cite{E}, G. Endimioni proved  that if $w$ be a word in a free
group such that finitely generated soluble groups in
$\mathcal{V}(w)$ are polycyclic, then every finitely generated
soluble group in  $\mathcal{V}(w^*)$ belongs to the variety
$\mathcal{V}(w)$. Before stating our next  result, we need  a
notation (see \cite{G1}). Let $\alpha$ be a non-zero  element of
some field of characteristic $p$. Denote the  group generated by
the matrices $\left\{ \begin{bmatrix} 1&0\\1&1\end{bmatrix},
 \begin{bmatrix} \alpha&0\\0&1\end{bmatrix}\right\}$ by $M(\alpha,p)$.\\
\par
{\bf Theorem D.} {\sl Let $w$
be a word in a free group  such
 that every infinitely presented $M(\alpha, p) \not\in \mathcal{V}(w)$ for all
$p\geq 0$ or $C_q\text{wr} C_{\infty}\not\in \mathcal{V}(w)$ for
all primes $q$. Then every infinite locally soluble group in
$\mathcal{V}(w^*)$ belongs to
the variety  $\mathcal{V}(w)$.}\\
\par
We note that
  the group  $M(\alpha,p)$ is finitely presented if and
only if \\
(i)~ $p\not = 0$ and $\alpha$ is algebraic over the prime field, or\\
(ii)~ $p=0$ and at least one of $\alpha$ or $\alpha^{-1}$ is an algebraic
integer  (see Lemma 11 of \cite{G1}).\\
\par
Theorem D generalizes Theorems 2 and 3  of \cite{E}, since we
note that if $\mathcal{V}$ is a variety  of groups in which every
finitely generated soluble group in $\mathcal{V}$ is polycyclic
then $\mathcal{V}$ contains no infinitely presented $M(\alpha, p)$
since $M(\alpha, p)$ is finitely generated metabelian;  the
subgroup $C_q^{(C_{\infty})}$ of $C_q \text{wr} C_{\infty}$ is not
finitely  generated and,  $C_q \text{wr} C_{\infty}$ is not
polycyclic for any prime $q$.\\
\par
\begin{center} {\bf\Large Proofs} \end{center}

We start the proof of  Theorem  A.\\
\par
{\em Proof of Theorem A.}
 Let $G$ be an infinite finitely generated locally graded restrained group in
$\mathcal{V}(w^*)$ and let $R$ be the finite residual of $G$.
Then $G/R$ is a finitely generated residually finite  group in
$\mathcal{V}(w^*)$ and so, by Lemma 1 of \cite{E},  it belongs to
$\mathcal{V}(w)$. Thus by hypothesis, $G/R$ is
polycyclic-by-finite. Therefore by repeated use of Lemma 3 of
\cite{KR}, $R$ is finitely generated. If $R$  is finite then $G$
is residually finite and so by Lemma 1 of \cite{E}, $G$ belongs
to the variety $\mathcal{V}(w)$. Now suppose, for a contradiction,
 that  $R$ is infinite. By hypothesis, $R$
has a  normal  proper subgroup of finite index   in $R$, then the
finite residual subgroup $T$  of $R$ is  proper in $R$. Therefore
$R/T$ is a residullay finite group in $\mathcal{V}(w)$ and so
$G/T$ is polycylic-by-finite. Thus $G/T$ is residually finite and
$R\subseteq T$, a contradiction. This completes the
proof. \;\;\;$\Box$\\
\par
 The following proposition  generalizes  the result of \cite{BP}.\\
\par
{\bf  Proposition 1.}
 {\sl Finitely generated residually finite groups in
a variety $\mathcal{V}$ in which every finite group  is nilpotent,
are nilpotent.\\}
\par
{\em Proof.}
We first prove that
 there exists a positive integer $k$ depending only on the variety $\mathcal{V}$
such that for all primes $p$, $C_p \text{wr} C_{p^k} \not\in
\mathcal{V}$.
 By the  Lemma of \cite{E1}, there exists an integer $t$ depending
only on $\mathcal{V}$ such that every $2$-generated metabelian
group in $\mathcal{V}$ is nilpotent of class at most $t$. Now
suppose that $C_p \text{wr} C_{p^m}  \in \mathcal{V}$ for some
prime $p$ and positive integer $m$. Since $C_p \text{wr} C_{p^m}$
is a $2$-generated
 metabelian group then it is nilpotent
 of class at most $t$. But the nilpotency class of
$C_p \text{wr} C_{p^m}$
 is exactly $p^m$, by a result of Liebeck (see \cite{L}
or Theorem 2.5 in page 76 of \cite{M}) and so $p^m\leq t$.
Now  the same argument as in     Theorem 2 of \cite{W} completes the proof.
 \;\;\;$\Box$\\
\par
Theorem A improves Theorem 3 of \cite{E} since
 by  the result  of \cite{BP}, in a variety, all finite
groups are nilpotent if and only if all finitely generated soluble
 groups are nilpotent. Therefore by Proposition 1, every  variety in which all
finitely generated soluble groups  are nilpotent  is contained
 in  a variety in
which all finitely generated residually finite groups are
polycyclic-by-finite.\\
\par
{\bf Corollary 2.} {\sl Let $w$ be a word in a free group such
that finitely generated soluble groups in $\mathcal{V}(w)$ are
nilpotent. Then every infinite locally graded restrained group in
 $\mathcal{V}(w^*)$ belongs to the  variety $\mathcal{V}(w)$.}\\
\par
{\em Proof.} As noticed before,  every finitely generated
residually finite group in $\mathcal{V}(w)$ is
polycyclic-by-finite. Let
 $G$ be an infinite locally graded restrained group in $\mathcal{V}(w^*)$ and
assume that
 $w$ is  a word in the free group of rank $n\in\mathbb{N}$. Let $x_1,
\dots, x_n \in G$, we must prove that $w(x_1, \dots, x_n)=1$.
Assume that there exists an infinite finitely generated subgroup
$H$ of $G$ which contains $x_1, \dots, x_n$. Then by Theorem A,
$H\in \mathcal{V}(w)$. Now, we may assume that every finitely
generated  subgroup of $G$  containing $x_1, \dots, x_n$ is
finite. Thus there exists an infinite locally finite subgroup $L$
which contains $x_1, \dots, x_n$ and so by Theorem 3 of \cite{E},
$L$ belongs to the variety $\mathcal{V}(w)$.
This completes the proof. \;\;\;$\Box$\\
\par
In the following lemmas we use some notion: we say that a word
$w\not=1$ in a free
group is a semigroup word
  if $w$ is of the form $uv^{-1}$, where $u$ and
$v$ are words in a free semigroup and we say, following \cite{LMR2}, that
 a group $G$ has no free
subsemigroups if and only if  for every pair (a,b) of elements of $G$,
the subsemigroup generated by $a, b$ has a relation of the form
$$(1)\;\;\;\;\;\;\; a^{r_1}b^{s_1}\cdots
 a^{r_j}b^{s_j}=b^{m_1}a^{n_1}\cdots b^{m_k}a^{n_k}$$
 where $r_i$, $s_i$, $m_i$ $n_i$ are all non-negative and $r_1$ and $m_1$ are
positive integers. If $(a,b)$ is a pair of elements in $G$ satsfying a
relation of type (1), then we call $j+k$ the width of the relation and the
 sum $r_1+\cdots +r_j+n_1+\cdots +n_k$
the exponent of $a$ (denoted  $\text{exp}(a)$) in the relation.\\
We say that a word  $w$ in a free group  $F$ generated by $x_1, \dots x_n$,
 is a commutator word  whenever $w$ belongs
to the derived subgroup of $F$. In the following we study
infinite groups in $\mathcal{V}(w^*)$ where $w$ is not a
commutator word. We note that if
 $w$ is not a commutator word then there  is a positive integer $e$
depending only on $w$ such that every group in the variety
$\mathcal{V}(w)$ is of exponent dividing $e$; for let $G$ be a
group in the variety generated by  a non-commutator word $w$,
since  $w$ is not a commuatator word,  for some $i$ the sum of the
exponents of $x_i$ in $w$ is non-zero: let this sum be $r$ and
let $g\in G$. If we replace $x_i$ by $g$ and $x_j$ by $1$ when
$j\not =i$, then $w$ assumes the value $g^r$. Thus $g$ has a
finite order $r$ and $G$ is of finite
exponent.\\
\par
{\bf Lemma 3.} {\sl Let $w$ be a semigroup word in the free group
of rank 2. Then every group in $\mathcal{V}(w^*)$ has no free
subsemigroups, and there exist  positive integers $M$ and $N$
depending only on $w$ such that for all pairs $(a,b)$ of elements
in $G$ there is a relation of the form
(1) whose width  and  $\text{exp}(a)$ is at most $M$ and $N$, respectively.}\\
\par
{\em Proof.} Let $a, b$ be in $G$. If $b$ is of finite order $m$ then
$ab^m=b^ma$,  $\text{exp}(a)=2$ and the width is $2$.
 Now,  assume that $b$ is of infinite
order and consider the two sets
  $X=\{a^{b^n} \;|\; n\in \mathbb{N}\}$ and $Y=\{b^m \;|\; m\in \mathbb{N}\}$.  If
$X$ is finite then the centre of $H=\left<a,b\right>$ is infinite
and so by Lemma 3 of \cite{E}, $H$ belongs to the variety
$\mathcal{V}(w)$. Therefore $w(a,b)=w(b,a)=1$ and so the pair
$(a,b)$ satisfies a relation of the form  (1) whose width and
 $\text{exp}(a)$ is
at most  $M_1$ and $N_1$, respectively, where $M_1$ and $N_1$
are  positive fixed integers
 depending only on $w$. Now we  may assume that
$X$ is infinite, then by the property $\mathcal{V}(w^*)$, there
exists a relation of the form
$$(a^{b^t})^{r_1}b^{s_1}\cdots
(a^{b^t})^{r_j}b^{s_j}=b^{m_1}(a^{b^t})^{n_1}\cdots
b^{m_k}(a^{b^t})^{n_k}$$
where  $r_i, s_i, m_i, n_i$ are  non-negative
 integers and $r_1, m_1, t$ are positive
integers; also
 the sum $r_1+\cdots +r_j+n_1+\cdots +n_k$ is the same $N_1$ and $j+k=M_1$.
 Therefore the pair
$(a,b)$ satisfies  a relation of the form (1) whose width is at most
$M:=\max\{2,M_1\}$ and  $\text{exp}(a)$ is at most
$N=\max\{2,N_1\}$. \;\;\;$\Box$\\
\par
Recall that a group $G$ is right orderable if there exists a total order
relation $\geq$ on $G$ such that for all $a, b, g$ in $G$, $a\geq b$ implies
$ag\geq bg$, equivalenty, if there exists a subset $P$ in $G$ such that
$PP=P$, $P\cup P^{-1}=G$, and $P\cap P^{-1}={1}$.\\
\par
{\bf Proposition 4.} {\sl Let $w$ be a semigroup word in the free
group of rank 2. Then every right orderable group in
$\mathcal{V}(w^*)$, belongs to the
variety $\mathcal{V}(w)$.}\\
\par
{\em Proof.} By  Theorem 5 of \cite{LMR2} and Lemma 3, $G$ is locally
nilpotent-by-finite. Let $x_1, \dots, x_n\in G$. Since $G$ is right
orderable, $G$ is torsion-free. Thus  every finitely generated subgroup of $G$
is  an infinite finitely generated nilpotent-by-finite group and so
residually finite. Therefore by Lemma 1 of \cite{E}, $G$ belongs to  the
variety $\mathcal{V}(w)$. \;\;\;$\Box$\\
\par
{\bf Lemma 5.} {\sl Let $w$ be a semigroup word in a free group. Then every
group in $\mathcal{V}(w^*)$ is restrained.}\\
\par
{\em Proof.} Let $G$ be a group in $\mathcal{V}(w^*)$ and let $x,
 y$ in $G$. We must
prove that $H=\left<x\right>^{\left<y\right>}$ is finitely generated. We may
assume that $y$ is of infinite order. Suppose that $w$ is in the free group of
rank $n>0$. Consider a partition of the set
 $X=\{xy^{-1}, xy^{-2},
\dots,\}$ in $n$ infinite subsets $X_1,X_2,\dots, X_n$.
 Then by the property $\mathcal{V}(w^*)$, there exist negative integers
$t_1, \dots, t_n$ such that
$$xy^{t_{f(1)}}\cdots xy^{t_{f(m)}}=xy^{t_{g(1)}}\cdots xy^{t_{g(s)}}$$
for some functions $f$ from $\{1,2,\dots, m\}$ to $\{1,2,\dots,n\}$ and $g$
from $\{1,2,\dots,s\}$ to $\{1,2,\dots,n\}$, where $m$ and $s$ depend
 only on $w$. Now,  arguing as in Lemma 1(ii) of
\cite{KR}, $H$ is finitely generated. This completes the proof. \;\;\;$\Box$\\
\par
{\bf Lemma 6.} {\sl  Let $w$ be a word in a free
 group such that $w$ is not a commutator word. Then every group in
$\mathcal{V}(w^*)$ is torsion. In particular, $G$ is restrained.}\\
\par
{\em Proof.} Let $G$ be a group. Suppose, for a contradiction,
that $G$  has an element $a$ of infinite order, then, by Lemma 3
of \cite{E}, $\left<a\right>$ belongs to the variety
$\mathcal{V}(w)$ and so $a$ is of finite
order, a contradiction. \;\;\;$\Box$\\
\par
We note that Theorem A
 can be applied for the following  words $w$ in a free
group: by Proposition 1 and the result of \cite{BP}, any word $w$ such that
 every finitely generated soluble group in the variety $\mathcal{V}(w)$ is
nilpotent; by Zelmanov's positive solution to the restricted Burnside problem
(see \cite{Z1} and \cite{Z2}), any non-commutator word $w$ and by Theorem A of
\cite{KR}, every semigroup word $w$.\\
By Theorem A and Lemmas 5 and  6 and  the above remarks we have\\
\par
{\bf Corollary 7.} {\sl Let $w$ be a non-commutator word or a
semigroup word in a free group. Then every infinite finitely
 generated locally graded group in $\mathcal{V}(w^*)$ belongs
to the variety $\mathcal{V}(w)$.}\\
\par
{\bf Lemma 8.} {\sl Let $G$ be an infinite  group in
$\mathcal{V}(w^*)$ and $H$ be a finite subgroup of $G$.
 If $G$ has an infinite normal locally soluble subgroup, then $H$
belongs to $\mathcal{V}(w)$.}\\
\par
{\em Proof.} Let $S$ be a normal locally soluble infinite subgroup of $G$.
If $S$ is $\Check{\text{C}}$ernikov, then $S$ has an  infinite normal
characteristic abelian subgroup (see \cite{R1} vol. I
 page 68) so $G$ has an infinite
normal abelian subgroup whence $G$ belongs to $\mathcal{V}(w)$ by
Lemma 3 of
\cite{E}.\\
Therefore, we may assume that $S$ is not $\Check{\text{C}}$ernikov.
By a result of Zaicev (see \cite{Z}), there is an infinite abelian subgroup
$B$ of $S$ such that $H$ normalizes $B$. Hence $B$ is an infinite normal
subgroup of the group $BH$ and so again by Lemma 3 of \cite{E},
 $H$ belongs to $\mathcal{V}(w)$.  \;\;\;$\Box$  \\
\par
{\em Proof of Theorem B.}
It suffices to
 prove that an infinite [(locally $\mathcal{P}$)-by-finite] group in
$\mathcal{V}(w^*)$ belongs to the variety $\mathcal{V}(w)$.
 Let $H$ be a normal locally $\mathcal{P}$-subgroup
 of $G$ of finite index. If $G$ is torsion, then $G$ is locally finite and $H$
is a locally soluble infinite normal subgroup of $G$, so by Lemma
8, $G \in \mathcal{V}(w)$. Therefore we may assume that $G$ has
an element $a$ of infinite order. Let $x_1, \dots, x_n$ be
arbitrary elements of $G$. Then $K=\left<a, x_1, \dots,
x_n\right>$ is a finitely generated $\mathcal{P}$-by-finite
infinite group and so by condition (3), $K \in \mathcal{V}(w)$. \;\;\;$\Box$\\
\par
{\bf Corollary 9.} {\sl Let $G$ be an infinite locally finite
$\mathcal{V}(w^*)$ group. If $G$ satisfies one of the following
conditions, then $G$ belongs to
the variety $\mathcal{V}(w)$.\\
(1)\; $G$ has an infinite locally soluble normal subgroup.\\
(2)\; $G$ contains an element with finite centralizer.\\
(3)\; $G$ contains an element of prime power order  with
 $\Check{\text{C}}$ernikov centralizer in $G$.\\}
\par
{\em Proof.}  Let $x_1, \dots, x_n$ be arbitrary elements of $G$, we must
prove that $w(x_1, \dots, x_n)=1$. Since $G$ is locally finite,
$H=\left<x_1, \dots, x_n\right>$ is finite.\\
If $G$ has an infinite locally soluble normal subgroup, then, by Theorem B,
$H \in \mathcal{V}(w)$.\\
If $G$ satisfies the conditions (2) or (3) then by Hartley's results of
\cite{H} and \cite{H2} $G$
is (locally soluble)-by-finite and so by part (1), the
proof is complete. \;\;\;$\Box$\\
\par
Let $w$ be a word in a free group. Now
 we state    some reductions in investigation of the equality $\mathcal{V}(w) \cup
\mathcal{F} =\mathcal{V}(w^*)$ on the class of  locally soluble
groups and locally
finite groups.\\
\par
Let $G$  be an infinite locally soluble
 group  in $\mathcal{V}(w^*)$. If $G$ is torsion
then by Corollary 9(1), $G$ belongs to the variety
$\mathcal{V}(w)$.
 Therefore we may assume
that $G$ has an element $g$  of infinite order and so in order to
prove that $G \in \mathcal{V}(w)$  it suffices to show that for
all $x_1, \dots, x_n$, the infinite finitely generated soluble
subgroup $\left<x_1, \dots, x_n, g\right>$
belongs to the variety $\mathcal{V}(w)$. Therefore we have\\
\par
{\bf Remark 10.} {\sl
Let $w$ be a  word
 in a free group. Then the following are equivalent:\\
(1)~ any infinite locally soluble group in $\mathcal{V}(w^*)$
belongs to the variety $\mathcal{V}(w)$.\\
(2)~  any infinite finitely generated soluble group  in
$\mathcal{V}(w^*)$  belongs to $\mathcal{V}(w)$.\\}
\par
We note that, by  Lemma 6, every finitely generated soluble group
in $\mathcal{V}(w^*)$, where $w$ is not a commutator word,
 is finite.\\
\par
Let $G$ be an infinite locally finite group in $\mathcal{V}(w^*)$.
 In order to prove
that $G \in \mathcal{V}(w)$,  we must show that $\left<x_1,
\dots, x_n\right> \in \mathcal{V}(w)$ for all $x_1, \dots, x_n \in
G$, therefore we may assume that $G$ is countable.  Fix $x_1,
\dots, x_n \in G$ and let $H=\left<x_1, \dots, x_n\right>$. If
$C_G(H)$ is infinite, then there is an infinite abelian subgroup
$A$ in $C_G(H)$, as $G$ is locally finite (see Theorem 3.43 of
\cite{R1}). therefore the centre of  $K=\left<A, H\right>$ is
infinite  and so by Lemma 3 of \cite{E},
 $K \in \mathcal{V}(w)$. Thus we may assume
that $C_G(H)$ is finite.  Also, by  Lemma 4 of \cite{E} and Corollary 9
 we may assume that $H$ is not
supersoluble and the centralizer of any
element in $G$ is infinite and the centralizer of every element of prime
power order is not $\Check{\text{C}}$ernikov. These conditions on a locally
finite group lead us to the following defenitions.\\
\par
We say that a group $G$ is an $\mathcal{L}$-group whenever
 $G$ is an infinite countable locally finite group and there exists a
 finite subgroup $H$ of $G$ such that\\
(1)\; $H$ is not supersoluble and $C_G(H)$ is finite.\\
(2)\; $C_G(x)$ is infinite for all $x \in G$.\\
(3)\; $C_G(g)$ is  not $\Check{\text{C}}$ernikov for all elements $g
\in G$ of prime power order.\\
(4)\; the largest normal locally soluble subgroup of $G$ is finite.\\
In this case, We say that $G$ is an $\mathcal{L}$-group with
respect to $H$. Also, we say that $G$ is an $\mathcal{L}^*$-group
with respect to $H$ whenever
 every infinite subgroup of $G$ which contains $H$, is an
 $\mathcal{L}$-group with respect to  $H$. By these discussions we have\\
\par
{\bf Remark 11.} {\sl Let $w$ be a word in a free group. Then the following
are equivalent:\\
1)~ an infinite locally finite group in
$\mathcal{V}(w^*)$, belongs to the variety $\mathcal{V}(w)$. \\
2)~  an infinite  $\mathcal{L}^*$-group in $\mathcal{V}(w^*)$,
belongs to the variety $\mathcal{V}(w)$.\\}
\par
We use Remark 11 for the study  of an infinite locally finite
group $G$ in $\mathcal{V}(w^*)$ where $w$ is not a commutator
word in the  free group of rank $n>0$, and
 obtain another condition on such  groups $G$. We prove that
$G$ is of finite exponent dividing $e$,
 where $e$ is a  positive integer depending only on $w$ such that
every group in the variety $\mathcal{V}(w)$ is of exponent
dividing $e$.
 For, let $a$ be an
element of $G$, then $C_G(a)$ is infinite and by Theorem 3.43
 of \cite{R1} there
exists an infinite abelian subgroup $A$ in $C_G(a)$. By Lemma 3
of \cite{E}, $A\in \mathcal{V}(w)$. Consider infinite subsets
$X_1=\dots=X_n=aA$. Therefore, by the property $\mathcal{V}(w^*)$,
there exist $a_1, \dots, a_n\in A$ such that $w(aa_1, \dots,
aa_n)=1$. Thus $w(a, \dots, a)w(a_1, \dots, a_n)=1$. But $w(a_1,
\dots, a_n)=1$ and so $w(a, \dots, a)=1$ and $a^e=1$.  Therefore
we
have:\\
\par
{\bf Remark 12.} {\sl Let $w$ be a non-commutator word
 in a free group and  $e$ be a  positive integer depending only on $w$ such
that every
 group in the variety $\mathcal{V}(w)$ is of exponent dividing
 $e$. Then the following
are equivalent:\\
1)~ any infinite locally finite group in
$\mathcal{V}(w^*)$ belongs to the variety $\mathcal{V}(w)$. \\
2)~  any infinite  $\mathcal{L}^*$-group of exponent dividing $e$
belongs to the variety $\mathcal{V}(w)$.}\\
\par
A natural question which arises  is the following: Is there an infinite
$\mathcal{L}^*$-group of finite exponent?\\
We only know that such a group is not simple. For by
 a result of L. G. Kov$\acute{\text{a}}$cs
 \cite{Ko}, any infinite, simple, locally
finite group $G$ involves infinitely many non-isomorphic non-abelian finite
simple groups; hence, if  $G$ satisfies non-trivial laws, then according to
a result of G. A. Jones (see Theorem of \cite{J}), the variety generated by
infinitely many finite simple groups is the variety of all groups. But the
variety generated by $G$ is a proper variety, a contradiction.\\
\par
Now we study infinite simple locally
 finite  groups in $\mathcal{V}(w^*)$, where
$w$ is a non-trivial word in a free group. As we have seen earlier,
 there is no infinite simple
locally finite  group  which satisfies  a non-trivial identity.
Call a simple locally finite group an $S$-group. The $S$-groups
 fall into two classes with widely different properties---the linear
 groups and the non-linear groups. Every linear $S$-group is a group of
 Lie type over an infinite locally finite field (see \cite{HS}).\\
\par
{\em Proof of Theorem C.} Suppose, for a contradiction, there exists
 a non-linear $S$-group $G$ in $\mathcal{V}(w^*)$. By
 a result of Hartley  \cite{H3},  there exists a section $C/D$ of $G$ such that
  $C/D$ is a direct product of finite alternating
 groups of unbounded orders. Thus $C/D$ is an infinite  residually finite
group in $\mathcal{V}(w^*)$ and so $C/D$ belongs to the variety
$\mathcal{V}(w)$. Since $C/D$ is  a direct product of finite
alternating groups of unbounded orders, the variety
$\mathcal{V}(w)$ contains infinitely many non-isomorphic finite
alternating groups. Therefore, by Theorem of \cite{J},
$\mathcal{V}(w)$ is the variety of all groups and so $w$ is  the
trivial word, a contradiction.
This completes the proof. \;\;\;\; $\Box$\\
\par
P. S. Kim in \cite{Kim} studied  $\mathcal{V}(w_2^*)$ on the
class of locally soluble groups, where
$w_2=[[x_1,x_2],[x_3,x_4]]$. For this word the
 variety $\mathcal{V}(w_2)$ is the  variety of metabelian groups. It is proved in
\cite{Kim}, that every  infinite locally soluble group in
$\mathcal{V}(w_2^*)$ is metabelian and also it is proved that any
infinite group belonging to $\mathcal{V}(w_2^*)$ is metabelian if
and only if there is no infinite simple group in
$\mathcal{V}(w_2^*)$. We study $\mathcal{V}(w^*)$ on the class of
locally finite groups, where $w$ is a soluble word that is
$w=w_d$ for some $d\in \mathbb{N}$ where $w_0=x$,
$w_i=[w_{i-1},w_{i-1}]$ and $w_{i-1}$
 is the  word on $2^{i-1}$ distinct letters which has been defined inductively,
 for all $i\in\mathbb{N}$.  \\
\par
{\bf Corollary 13.} {\sl Let $w$ be a soluble word and let   $G$
be an infinite locally finite $\mathcal{V}(w^*)$-group.
 Then  the following are equivalent:\\
(1)~ $G\in\mathcal{V}(w)$.\\
(2)~ $G$ has no infinite linear simple locally finite section.\\}
\par
{\em Proof.} Suppose that (1) is true. Then $G$ is soluble and
(2) is clear. Now suppose that (2) is true and $w=w_d$ for some
positive integer $d$. Suppose, for a contradiction, that
$G\not\in\mathcal{V}(w)$. Thus $G$ is not soluble of derived
length at most $d$. Suppose, if possible, that $K=G^{(d+1)}$ is
finite. Then $H=G^{(d)}$  is an FC-group and so $H$ is soluble by
applying suitably Lemma 1 of \cite{AT2}. Thus $G$ is a torsion
soluble group and so by Theorem B,  $G$ is soluble of derived
length at most $d$,
 a contradiction. Hence
$K$ is infinite and so $G/K$ is a soluble group of derived length $d$.
Therefore $G^{(d)}=G^{(d+1)}$ that is $H=H'$, which implies that $H$ is a
perfect group. Suppose  that $H$ has an infinite  proper normal subgroup $N$,
then $H/N$ is soluble of derived length at most
 $d$, this implies $H=H^{(d)}\leq N$ since $H$ is perfect, a
contradiction. Let $N$
be a finite normal subgroup of
$H$, then $C_H(N)$ has finite index in $H$. Since $H$ has no infinite normal
proper subgroups, $C_H(N)=H$. Hence the centre $Z$ of $H$
 is the unique maximal
normal subgroup of $H$ so that $S=H/Z$
 is simple. By Theorem C, $S$ is an infinite linear simple locally
finite group, which is a contradiction. \;\;\;\; $\Box$\\
\par
Now, we start proving Theorem D, for this  we need the following lemma:\\
\par
{\bf Lemma 14.} {\sl Every infinite locally soluble group of finite rank in
 $\mathcal{V}(w^*)$  belongs to the variety $\mathcal{V}(w)$.}\\
\par
{\em Proof.} Let $G$ be an infinite locally soluble group of
finite rank in $\mathcal{V}(w^*)$. By
 Remark 10, we may assume that $G$ is finitely generated.
Therefore $G$   is a minimax group, and so by Theorem 10.33 of
\cite{R1},  the finite residual $R$ of $G$ is the direct product of finitely
many quasicyclic subgroups of $G$, thus $G$ is residually finite or $G$ has an
infinite normal abelian subgroup, then, by Lemma 1 or Lemma 3 of \cite{E}
respectively,  the proof is complete.\\
\par
{\em Proof of Theorem D.} Let $G$ be an infinite locally soluble
group in $\mathcal{V}(w^*)$. By Remark 10,
 we may assume that $G$ is a finitely generated infinite soluble
group. Firstly, suppose  that $w$ is a word such that
$C_p\text{wr} C_{\infty}\not\in \mathcal{V}(w)$ for all primes
$p$. We prove that $G$ is a minimax group and so $G$ is of finite
rank,
then  Lemma 14 completes     the proof.\\
By a  deep result of Kropholler (see \cite{K}), which asserts
that every finitely generated soluble group having no sections of
type $C_p \text{wr} C_{\infty}$ is minimax, it suffices to show
that if $C_p \text{wr} C_{\infty} \in \mathcal{V}(w^*)$ then $C_p
\text{wr} C_{\infty} \in \mathcal{V}(w)$. But $C_p \text{wr}
C_{\infty}$  has an infinite normal abelian subgroup, therefore
by Lemma 3 of \cite{E}, $C_p \text{wr} C_{\infty} \in
\mathcal{V}(w)$, which  is a
contradiction.\\
Now, suppose that $w$ is a word such that
  every infinitely presented $M(\alpha, p) \not\in \mathcal{V}(w)$.
If $G$ is not
semi-polycyclic group (see \cite{G1}),
 then there exists a subgroup $H$ of a quotient group of
$G$ which is isomorphic to an infinitely  presented  $M(\alpha,
p)$. But $M(\alpha, p)$  is an infinite residually finite group
in $\mathcal{V}(w^*)$
 and so $M(\alpha, p) \in \mathcal{V}(w)$, a contradiction. Therefore $G$ is
semi-polycyclic and so is of finite rank (see   \cite{G1}). Thus,
by Lemma 14,  $G\in\mathcal{V}(w)$. This completes the proof.
\;\;\;\;
$\Box$\\


\begin{thebibliography}{99}
\bibitem{A1} A. Abdollahi, {\sl A characterization of infinite 3-abelian
groups,} Arch. Math. (Basel), {\bf 73} (1999), 104-108.
\bibitem{A2} A. Abdollahi, {\sl Finitely generated soluble groups with an
Engel condition on infinite subsets,}  Rend. Sem. Mat. Univ.
Padova {\bf 103} (2000) 47-49.
\bibitem{A3} A. Abdollahi, {\sl Some Engel conditions on infinite subsets of
certain groups,}  Bull. Austral. Math. Soc. {\bf 62} (2000) 141-148.
\bibitem{AT2} A. Abdollahi and B. Taeri, {\sl A condition on a certain variety
of groups,}  Rend. Sem. Mat. Univ. Padova {\bf 104} (2000), 129-134.
\bibitem{AT1} A. Abdollahi and B. Taeri, {\sl A condition on finitely
generated soluble groups,}  Comm. Algebra {\bf 27} (1999), 5633-5638.
\bibitem{AT3} A. Abdollahi and B. Taeri, {\sl Some conditions on infinite
subsets of infinite groups,}  Bull. Malaysian Math. Soc. (2) {\bf 22} no. 1
(1999) 87-93.
\bibitem{BP} M. Boffa and F. Point, {\sl Identit$\acute{\text{e}}$s de Engel
g$\acute{\text{e}}$n$\acute{\text{e}}$ralis$\acute{\text{e}}$es,}
 C. R. Acad. Sci. Paris, {\bf 313} (1991), 909-911.
\bibitem{Bouk} A. Boukaroura, {\sl A condition of infinite groups for
satisfying certain laws,} to appear in Algebra Colloq.
\bibitem{D1} C. Delizia, {\sl Finitely generated soluble groups with a
condition on infinite subsets,} Istit. Lombardo Accad. Sci. Lett. Rend. A {\bf
128} (1994), 201-208.
\bibitem{D2} C. Delizia, {\sl On groups with a nilpotence condition on
infinite subsets,} Algebra Colloq. {\bf 2} (1995), 97-104.
\bibitem{D3} C. Delizia, {\sl On certain residually finite groups,}
Comm.Algebra {\bf 24} (1996), 3531-3535.
\bibitem{DRS} C. Delizia, A.  Rhemtulla and H. Smith, {\sl Locally graded
groups with a nilpotency condition on infinite subsets,} to appear.
\bibitem{E1} G. Endimioni, {\sl Conditions de nilpotence dans certaines
vari$\acute{e}$t$\acute{e}$s de groupes,} C. R. Acad. Sci. Paris, {\bf 310}
(1990), 325-327.
\bibitem{E}  G. Endimioni, {\sl
 On a combinatorial problem in varieties of groups,} Comm. Algebra {\bf 23}
(1995), 5297-5307.
\bibitem{G} J. R. J. Groves, {\sl A conjecture of Lennox and
Wiegold concerning
supersoluble groups,} J. Austral. Math. Soc. (Series A) {\bf 35}
(1983),218-220.
\bibitem{G1} J. R. J. Groves, {\sl Soluble groups in which every finitely
generated subgroup is finitely presented,} J. Austral. Math. Soc. (Series A)
{\bf 26} (1978), 115-125.
\bibitem{Gru} K. W. Gruenberg, {\sl Two theorems on Engel groups,} Proc.
Cambridge Philos. Soc. {\bf 49} (1953), 377-380.
\bibitem{J} G. A. Jones, {\sl Varieties and simple groups,} J. Austral. Math.
Soc. {\bf 17} (1974), 163-173.
\bibitem{Kim} P. S. Kim, {\sl A condition for locally soluble groups to be
metabelian,} Houston J. Math. {\bf 20} (1994), 193-199.
\bibitem{KR} Y. K. Kim and A.  Rhemtulla, {\sl Weak maximality condition and
polycyclic groups,} Proc. Amer. Math. Soc. {\bf 123} (1995), 711-714.
\bibitem{KRS} P. S. Kim, A.  Rhemtulla and H. Smith, {\sl A characterization
of infinite metabelian groups,} Houston J. Math. {\bf 17} (1991), 429-437.
\bibitem{Ko} L. G. Kov$\acute{\text{a}}$cs,
 {\sl Varieties and finite groups,} J.
Austral. Math. Soc.  {\bf 10} (1969), 5-19.
\bibitem{K}  P. H. Kropholler, {\sl On
 finitely generated soluble groups with no large wreath product sections,}
Proc. London Math. Soc. {\bf 49} (1984), 155-169.
\bibitem{LW} J. C. Lennox and J. Wiegold, {\sl Extensions of a problem of Paul
Erd\"os on groups,} J. Austral. Math. Soc. {\bf 31} (1981), 459-463.
\bibitem{L} H. Liebeck, {\sl Concerning nilpotent wreath products,} Proc.
Cambridge Philos. Soc., {\bf 58} (1962), 443-451.
\bibitem{LM} P. Longobardi and M. Maj, {\sl A finiteness condition concerning
commutators in groups,} Houston J. Math. {\bf 19} (1993), 505-512.
\bibitem{LM2} P. Longobardi and M. Maj, {\sl Finitely generated soluble groups
with an Engel condition on infinite subsets,} Rend. Sem. Mat. Univ. Padove
{\bf 89} (1993), 97-102.
\bibitem{LMMR} P. Longobardi, M. Maj, A. Mann and A.  Rhemtulla, {\sl Groups
with many nilpotent subgroups,} Rend. Sem. Mat. Univ. Padova {\bf 95} (1996),
143-152.
\bibitem{LMR} P. Longobardi, M. Maj and A.  Rhemtulla, {\sl Infinite groups
in a given variety and Ramsey's theorem,} Comm. Algebra {\bf 20} (1992),
127-139.
\bibitem{LMR2} P. Longobardi, M. Maj and A. Rhemtulla, {\sl Groups with no
free subsemigroups,} Trans. Amer. Math. Soc. {\bf 347} (1995), 1419-1427.
\bibitem{H2}  B. Hartley, {\sl Fixed points of automorphisms of prime power
order of locally finite groups and Chevalley groups,} J. London
Math. Soc., {\bf 37} (1988), 421-436.
\bibitem{H3} B. Hartley, {\sl  Centralizing properties in simple locally finite
groups and large finite  classical groups,}
 J. Austral. Math. Soc. (Series A) {\bf 49} (1990), No. 3, 502-513.
\bibitem{H}  B. Hartley, {\sl A general Brauer-Fowler Theorem and centralizers
in locally finite groups,} Pacific J. Math. {\bf 152} (1992), 101-117.
\bibitem{HS} B. Hartley and G. Shute, {\sl Monomorphisms and direct limits of
finite groups of Lie type,} Quart. J. Math. Oxford (2) {\bf 35} (1984), 49-71.
\bibitem{Mac} I. D. Macdonald, {\sl On certain varieties of groups II,}  Math.
Zeitschr. {\bf 78} (1962), 175-188.
\bibitem{M} J. D. P. Meldrum, {\sl Wreath products of groups and semigroups,}
Pitman Monographs and Surveys in Pure and applied Mathematics, {\bf 74}
London, 1995.
\bibitem{N} B. H. Neumann, {\sl A problem of Paul Erd\"os on groups,} J.
Austral. Math. Soc. (Series A) {\bf 21} (1976), 467-472.
\bibitem{PS} O. Puglisi and L. S. Spiezia, {\sl A combinatorial property on
certain infinite groups,}  Comm. Algebra {\bf 22} (1994), 1457-1465.
\bibitem{RH} A. Rhemtulla and H. Smith, {\sl On infinite locally finite
groups,} Canad. Math. Bull. {\bf 37} (1994), 537-544.
\bibitem{R1}  D. J. S. Robinson, {\sl Finiteness conditions and generalized
soluble groups, I, II,} Springer-Verlag, Berlin, 1972.
\bibitem{S1} L. S. Spiezia, {\sl A property of the variety of 2-Engel groups,}
Rend. Sem. Mat. Uinv. Padova {\bf 91} (1994), 225-228.
\bibitem{S2} L. S. Spiezia, {\sl  A characterization of third Engel groups,}
Arch. Math. (Basel), {\bf 64} (1995), 369-373.
\bibitem{Taeri} B. Taeri, {\sl A combinatorial condition on a ceratin
variety of groups,} Arch. Math. (Basel)
\bibitem{W} J. S. Wilson, {\sl Two-generator conditions for residually finite
groups,} Bull. London Math. Soc. {\bf 23} (1991), 239-248.
\bibitem{Z} D. I. Zaicev, {\em On Solvable Subgroups of Locally Solvable
Groups, } Dokl. Akad. Nauk SSSR {\bf 214} (1974), 1250-1253, translation in
Soviet Math. Dokl. 15 (1974), 342-345.
\bibitem{Z1} E. I. Zelmanov, {\sl The solution of the restricted Burnside
problem for groups of odd exponent,} Math. USSR Izv. {\bf 36} (1991), 41-60.
\bibitem{Z2} E. I. Zelmanov, {\sl The solution of the restricted Burnside
problem for  2-groups,} Math. Sb.  {\bf 182} (1991), 568-592.

\end{thebibliography}
\end{document}